\input amstex
\documentstyle{amsppt}
\topmatter \magnification=\magstep1 \pagewidth{5.2 in}
\pageheight{6.7 in}
\abovedisplayskip=10pt \belowdisplayskip=10pt
\parskip=8pt
\parindent=5mm
\baselineskip=2pt
\title
 On  $p$-adic $q$-$l$-functions and sums of powers
\endtitle
\author   Taekyun Kim    \endauthor

\affil{ {\it Jangjeon Research Institute for Mathematical Sciences $\&$ Physics,\\
        Ju-Kong Building 103-Dong 1001-ho,\\
        544-4 Young-chang Ri Hapcheon-Up Hapcheon-Gun Kyungnam,\\
         678-802, S. Korea\\
        e-mail: tkim64$\@$hanmail.net ( or tkim$\@$kongju.ac.kr)}}\endaffil
\thanks  2000 Mathematics Subject Classification:  11S80, 11B68, 11M99 .\endthanks
\abstract{ In this paper, we give an explicit $p$-adic expansion
of $$\sum_{\Sb j=1\\ (j,p)=1\endSb}^{np} \frac{(-1)^j }{[j]_q^r}
$$  as a power series in $n$. The coefficients are values of
$p$-adic $q$-$l$-function for $q$-Euler numbers.
 }\endabstract
\rightheadtext{  On the $p$-adic interpolation function}
\leftheadtext{T. Kim}
\endtopmatter

\document

\head \S 1. Introduction \endhead Let $p$ be a fixed prime.
Throughout this paper  $\Bbb Z_p,\,\Bbb Q_p , \,\Bbb C$ and $\Bbb
C_p$ will, respectively, denote the ring of $p$-adic rational
integers, the field of $p$-adic rational numbers, the complex
number field and the completion of algebraic closure of $\Bbb Q_p
, $ cf.[1, 4, 6, 10].  Let $v_p$ be the normalized exponential
valuation of $\Bbb C_p$ with $|p|_p=p^{-v_p(p)}=p^{-1}.$ When one
talks of $q$-extension, $q$ is variously considered as an
indeterminate, a complex number $q\in\Bbb C ,$ or a $p$-adic
number $q\in\Bbb C_p$. If $q\in\Bbb C ,$ one normally assumes
$|q|<1 .$  If $ q \in \Bbb C_p ,$  then we assume $|q-1|_p <
p^{-\frac1{p-1}},$ so that $q^x=\exp(x\log q)$ for $|x|_p \leq 1.$
Kubota and Leopoldt proved the existence of meromorphic functions,
$L_p(s, \chi )$, defined over the $p$-adic number field, that
serve as $p$-adic equivalents of the Dirichlet $L$-series, cf.[10,
11]. These $p$-adic $L$-functions interpolate the values
$$L_p(1-n, \chi)=-\frac{1}{n}(1-\chi_n(p)p^{n-1})B_{n,\chi_n},
\text{ for $n\in\Bbb N=\{1, 2,\cdots, \}  $ ,}$$ where
$B_{n,\chi}$ denote the $n$th generalized Bernoulli numbers
associated with the primitive Dirichlet character $\chi ,$ and
$\chi_n=\chi w^{-n} ,$ with $w$  the $Teichm\ddot{u}ller$
character, cf.[8, 10]. In [10], L. C. Washington have proved the
below interesting formula:
$$\sum_{\Sb
j=1\\(j,p)=1\endSb}^{np}\frac{1}{j^r}=-\sum_{k=1}^{\infty}\binom{-r}{k}(pn)^kL_p(r+k,
w^{1-k-r}), \text{ where $\binom{-r}{k} $ is binomial
coefficient}.
$$ To give the $q$-extension of the above Washington result,
author derived the sums of powers of consecutive $q$-integers as
follows:
$$\sum_{l=0}^{n-1}q^l[l]_q^{m-1}=\frac{1}{m}\sum_{l=0}^{m-1}\binom{m}{l}q^{ml}\beta_l[n]_q^{m-l}
+\frac{1}{m}(q^{mn}-1)\beta_m, \text{ see [6, 7] ,} \tag*$$ where
$\beta_m$ are $q$-Bernoulli numbers. By using (*), we gave an
explicit $p$-adic expansion
$$\aligned
&\sum_{\Sb j=1\\(j,p)=1 \endSb}^{np}
\frac{q^j}{[j]_q^r}=-\sum_{k=1}^{\infty}\binom{-r}{k}[pn]_q^k
L_{p,q}(r+k, w^{1-r-k})\\
&-(q-1) \sum_{k=1}^{\infty}\binom{-r}{k}[pn]_q^k T_{p,q}(r+k,
w^{1-r-k})-(q-1)\sum_{a=1}^{p-1}B_{p,q}^{(n)}(r,a:F),
\endaligned$$
where $L_{p,q}(s, \chi)$ is $p$-adic $q$-$L$-function (see [7] ).
Indeed, this is a $q$-extension result due to Washington,
corresponding to the case $q=1$, see [10]. For a fixed positive
integer $d$ with $(p,d)=1$, set
$$\split &
X =X_d =\varprojlim_N \Bbb Z/dp^N ,\cr &  X_1 = \Bbb Z_p , X^\ast
=\bigcup_{\Sb  0<a <dp\\(a,p)=1
\endSb} a +d p\Bbb Z_p , \\
& a +d p^N \Bbb Z_p = \{ x\in X | x \equiv a
\pmod{p^N}\},\endsplit$$ where $a\in\Bbb Z$ satisfies the
condition $0\leq a < d p^N$, (cf.[3, 4, 9]). We say that $f$ is a
uniformly differentiable function at a point $a\in\Bbb Z_p$, and
write $f\in UD(\Bbb Z_p )$, if the difference quotients $F_f (x,y)
= \frac{f(x) -f(y)}{x-y}$ have a limit $f^\prime (a)$ as $(x,y)\to
(a,a)$, cf.[3]. For $f\in UD (\Bbb Z_p )$, let us begin with the
expression
$$
\frac{1}{[p^N ]_q } \sum_{0\leq j< p^N } q^j f(j) =\sum_{0\leq j<
p^N } f(j) \mu_q (j+ p^N \Bbb Z_p ), \text{ cf.[1, 3, 4, 7, 8,
9],}
$$
which represents  a $q$-analogue of Riemann sums for $f$. The
integral of $f$ on $\Bbb Z_p$ is defined as the limit of those
sums(as $n\to \infty$) if this limit exists. The $q$-Volkenborn
integral of a function $f \in UD(\Bbb Z_p )$ is defined by
$$I_q (f) = \int_X f(x) d\mu_q (x) =\int_{X_d} f(x) d\mu_q (x)=
\lim_{N\to\infty} \frac{1}{[dp^N]_q} \sum_{x=0}^{dp^N -1}f(x) q^x
, \text{ cf. [2, 3] }.\tag 1$$ It is well known that the familiar
Euler polynomials $E_n(z)$ are defined by means of the following
generating function:
$$F(z,
t)=\frac{2}{e^t+1}e^{zt}=\sum_{n=0}^{\infty}E_n(z)\frac{t^n}{n!},
\text{ cf.[1, 5].}$$ We note that, by substituting $z=0$,
$E_n(0)=E_n$ are the familiar $n$-th Euler numbers. Over five
decades ago, Carlitz defined $q$-extension of Euler numbers and
polynomials, cf.[1, 4, 5]. Recently, author  gave another
construction of $q$-Euler numbers and polynomials (see [1, 5, 9]).
By using author's $q$-Euler numbers and polynomials, we gave the
alternating sums of powers of consecutive $q$-integers as follows:
$$2\sum_{l=0}^{n-1}(-1)^l[l]_q^m=(-1)^{n+1}\sum_{l=0}^{m-1}
\binom ml
q^{nl}E_{l,q}[n]_q^{m-l}+\left((-1)^{n+1}q^{nm}+1\right)E_{m,q},
$$ where $E_{l,q}$ are $q$-Euler numbers (see [5] ). From this
result, we can study the $p$-adic interpolating function for
$q$-Euler numbers  and sums of powers due to author [7].
Throughout this paper, we use the below notation:
$$\aligned
&[x]_q=\frac{1-q^x}{1-q}=1+q+q^2+\cdots +q^{x-1},\\
&[x]_{-q}=\frac{1-(-q)^x}{1-q}=1-q+q^2-q^3+\cdots+(-q)^{x-1},
\text{ cf.[5, 9]. }
\endaligned$$
Note that when $p$ is prime $[p]_q$ is an irreducible polynomial
in $Q[q].$ Furthermore, this means that $Q[q]/[p]_q$ is a field
and consequently rational functions $r(q)/s(q)$ are well defined
$\mod [p]_q$ if $(r(q), s(q))=1 $.
 In a recent paper [5] the author
constructed the new $q$-extensions of Euler numbers and
polynomials.
 In Section 2, we introduce the $q$-extension of Euler numbers and polynomials.
In Section 3 we construct a new $q$-extension of Dirichlet's type
$l$-function which interpolates the $q$-extension of generalized
Euler numbers attached to $\chi$ at negative integers.
 The values of this function at negative
integers are algebraic, hence may be regarded as lying in an
extension of $\Bbb Q_p$. We therefore look for a $p$-adic function
which agrees with at negative integers. The purpose of this paper
is to construct the new $q$-extension of generalized Euler numbers
attached to $\chi$ due to author and prove the existence of a
specific $p$-adic interpolating function which interpolate the
$q$-extension of generalized Bernoulli polynomials at negative
integer. Finally, we give an explicit $p$-adic expansion
$$\sum_{\Sb j=1 \\ (j, p)=1 \endSb}^{np}\frac{(-1)^j}{[j]_q^r
}, $$ as a power series in $n$. The coefficients are values of
$p$-adic $q$-$l$-function for $q$-Euler numbers.

 \head 2. Preliminaries  \endhead

For any non-negative integer $m$, the $q$-Euler numbers,
$E_{m,q}$, were represented by
$$\frac{2}{[2]_q}\int_{\Bbb Z_p}q^{-x}[x]_q^m d\mu_{-q}(x)=E_{m,q}=2\left(\frac{1}{1-q}\right)^m\sum_{i=0}^m\binom mi
(-1)^i\frac{1}{1+q^{i}}, \text{  see [9] }. \tag2$$ Note that
$\lim_{q\rightarrow 1} E_{m,q}=E_m .$ From Eq.(2), we can derive
the below generating function:
$$F_q(t)=2
e^{\frac{t}{1-q}}\sum_{j=0}^{\infty}\frac{1}{1+q^{j}}(-1)^j(\frac{1}{1-q})^j\frac{t^j}{j!}
=\sum_{j=0}^{\infty}E_{n,q} \frac{t^n}{n!}. \tag 3$$ By using
$p$-adic $q$-integral, we can also consider the $q$-Euler
polynomials, $E_{n, q}(x)$, as follows:
$$E_{n,q}(x)=\frac{2}{[2]_q}\int_{\Bbb Z_p}q^{-t}[x+t]_q^n d\mu_{-q}(t)
=2(\frac{1}{1-q})^n\sum_{k=0}^n\binom nk \frac{(-q^x)^k}{1+q^{k}},
\text{ see [5, 9].}\tag4$$ Note that
$$E_{n,q}(x)=\frac{2}{[2]_q}\int_{\Bbb Z_p}([x]_q+q^x[t]_q)^n q^{-t}d\mu_{-q}(x)
=\sum_{j=0}^n\binom nj q^{jx}E_{j,q}[x]_q^{n-j}. \tag5$$ By (4),
we easily see that
$$\sum_{n=0}^{\infty}E_{n,q}(x)\frac{t^n}{n!}=F_q(x,t)=2e^{\frac{t}{1-q}}
\sum_{j=0}^{\infty}\frac{(-1)^j}{1+q^{j}}q^{jx}(\frac{1}{1-q})^j\frac{t^j}{j!}.
\tag6$$ From (6), we derive
$$F_q(x, t)=2\sum_{n=0}^{\infty}(-1)^ne^{[n+x]_q
t}=\sum_{n=0}^{\infty}E_{n,q}(x)\frac{t^n}{n!}. \tag7$$

 \head 3. On the $q$-analogue of Hurwitz's type $\zeta$-function associated with
$q$-Euler numbers \endhead

In this section, we assume that $q\in\Bbb C$ with $|q|<1$. It is
easy to see that
$$E_{n,q}(x)=[m]_q^n\sum_{a=0}^{m-1}(-1)^a
E_{n,q^m}(\frac{a+x}{m}), \text{ see [1, 5] }, $$ where $m$ is odd
positive integer. From (7), we can easily derive the below
formula:
$$E_{k,q}(x)=\frac{d^k}{dt^k}F_q(x,t) |_{t=0}=2\sum_{n=0}^{\infty}(-1)^n[n+x]_q^k. \tag8$$
Thus, we can consider a $q$-$\zeta$-function which interpolates
$q$-Euler numbers at negative integer as follows:
  \proclaim{Definition 1} For $s\in\Bbb C$, define
  $$\zeta_{E,q}(s,
  x)=[2]_q\sum_{m=1}^{\infty}\frac{(-1)^n}{[n+x]_q^s}.$$
Note that $\zeta_{E,q}(s,x)$ is meromorphic function in whole
complex plane.
\endproclaim
By using Definition 1 and Eq.(8), we obtain the following:
\proclaim{Proposition 2} For any positive integer $k$, we have
$$\zeta_{E,q}(-k, x)=E_{k,q}(x).$$
\endproclaim
Let $\chi$ be the Dirichlet character with conductor $f\in\Bbb N$.
Then we define the generalized $q$-Euler numbers attached to
$\chi$ as
$$F_{q,\chi}(t)=2\sum_{n=0}^{\infty}e^{[n]_qt}\chi(n)(-1)^n=\sum_{n=0}^{\infty}
E_{n,\chi, q}\frac{t^n}{n!}. \tag9$$ Note that
$$E_{n,\chi,q}=[f]_q^n\sum_{a=0}^{f-1}\chi(a)(-1)^aE_{n,q^f}(\frac{a}{f}),
\text{ where $f(=odd) \in\Bbb N$ }.\tag10$$ By (9),  we easily see
that
$$\frac{d^k}{dt^k}F_{q,\chi}(t)|_{t=0}=E_{k,\chi,q}=2\sum_{n=1}^{\infty}\chi(n)(-1)^n[n]_q^k
\tag11 $$

\proclaim{Definition 3} For $s\in\Bbb C$, we define Dirichlet's
type $l$-function as follows:
$$l_q(s,
\chi)=[2]_q\sum_{n=1}^{\infty}\frac{\chi(n)(-1)^n}{[n]_q^s}. $$
\endproclaim
From (11) and Definition 3, we can derive the below theorem.
\proclaim{ Theorem 4} For $k \geq 1$, we have
$$l_q(-k, \chi)=E_{k,\chi,q}. $$
\endproclaim
In [5], it was known that
$$2\sum_{l=0}^{n-1}(-1)^l[l]_q^m=\left((-1)^{n+1}q^nE_{m,q}(n)+E_{m,q}
\right), \text{ where $m, n \in\Bbb N.$} \tag12 $$ From (4) and
(12), we derive
$$\aligned
&2\sum_{l=0}^{n-1}(-1)^l[l]_q^m\\
&=(-1)^{n+1}\sum_{l=0}^{m-1}\binom ml
q^{nl}E_{l,q}[n]_q^{m-l}+\left((-1)^{n+1}q^{nm}+1)\right)E_{m,q}.\endaligned\tag13$$
Let $s$ be a complex variable, and let $a$ and $F(=odd)$ be the
integers with $0<a<F.$ We now consider the partial $q$-zeta
function as follows:
$$H_q(s,a:F)=\sum_{\Sb m\equiv a(F)\\ m>0
\endSb}\frac{(-1)^m}{[m]_q^s}=(-1)^a\frac{[F]_q^{-s}}{2}\zeta_{E, q^F}(s, \frac{a}{F}).\tag14 $$
For $n\in\Bbb N$, we note that
$H_q(-n,a:F)=(-1)^a\frac{[F]_q^n}{2}E_{n,q^F}(\frac{a}{F}).$ Let
$\chi$ be the Dirichlet's character with conductor $F(=odd)$. Then
we have
$$l_q(s,\chi)=2\sum_{a=1}^F\chi(a)H_q(s,a:F). \tag 15$$
The function $H_q(s,a:F)$ will be called the $q$-extension of
partial zeta function which interpolates $q$-Euler polynomials at
negative integers. The values of $l_q(s, \chi)$ at negative
integers are algebraic, hence may be regarded as lying in an
extension of $\Bbb Q_p$. We therefore look for a $p$-adic function
which agrees with $l_q(s,\chi)$ at the negative integers in
Section 4.

\head \S 4. $p$-adic $q$-$l$-functions and sums of powers \endhead

We define $<x>=<x:q>=\frac{[x]_q}{w(x)}$, where $w(x)$ is the
$Teichm\ddot{u}ller$ character. When $F(=odd)$ is multiple of $p$
and $(a,p)=1$, we define a $p$-adic analogue of (14) as follows:
$$H_{p,q}(s,
a:F)=\frac{(-1)^a}{2}<a>^{-s}\sum_{j=0}^{\infty}\binom{-s}{j}q^{ja}
\left(\frac{[F]_q}{[a]_q}\right)^j E_{j,q^{F}}, \text{ for
$s\in\Bbb Z_p$.} \tag16$$ Thus, we note that
$$\aligned
H_{p,q}(-n, a:F)&=\frac{(-1)^a}{2}<a>^n\sum_{j=0}^n\binom
nj q^{ja}\left(\frac{[F]_q}{[a]_q}\right)^jE_{j,q^F}\\
&=\frac{(-1)^a}{2}w^{-n}(a)[F]_q^nE_{n,q^F}(\frac{a}{F})=w^{-n}(a)H_q(-n,a:F),
\text{ for $n\in\Bbb N $. }
\endaligned \tag17$$
We now construct the $p$-adic analytic function which interpolates
$q$-Euler number at negative integer as follows:
$$l_{p,q}(s, \chi)=2\sum_{\Sb a=1\\ (a,p)=1
\endSb}^F\chi(a)H_{p,q}(s, a:F). \tag18$$
In [5, 9], it was known that
$$E_{k,\chi,q}=\frac{2}{[2]_{q^f}}\int_{X}\chi(x)[x]_q^kq^{-x} d\mu_{-q}(x), \text{ for
$k\in\Bbb N $.}$$ For $f(=odd)\in\Bbb N,$ we note that
$$E_{n,\chi,q}=[f]_q^n\sum_{a=0}^{f-1}\chi(a)(-1)^aE_{n,q^f}(\frac{a}{f}).$$

Thus, we have
$$\aligned
l_{p,q}(-n, \chi)&=2\sum_{\Sb a=1 \\ (p,a)=1
\endSb}^F\chi(a)H_{p,q}(-n,a:F)=\frac{2}{[2]_{q^f}}\int_{X^*}\chi w^{-n}(x)[x]_q^nq^{-x}d\mu_{-q}(x)\\
&=E_{n,\chi w^{-n},q}-[p]_q^n\chi w^{-n}(p)E_{n, \chi w^{-n}, q^p
}.\endaligned\tag 18-1$$ In fact,
$$l_{p,q}(s,\chi)=2\sum_{a=1}^F(-1)^a<a>^{-s}\chi(a)\sum_{j=0}^{\infty}
\binom{-s}{j}q^{ja}\left(\frac{[F]_q}{[a]_q}\right)^jE_{j,q^F},
\text{ for $s\in\Bbb Z_p$.}\tag19$$ This is a $p$-adic analytic
function and has the following properties for $\chi=w^t$:
$$l_{p,q}(-n, w^t)=E_{n,q}-[p]_q^nE_{n, q^p}, \text{ where $n
\equiv t$ ($\mod p-1$)}, \tag20$$
$$l_{p,q}(s,t)\in\Bbb Z_p \text{ for
all $s\in\Bbb Z_p$ when $t\equiv 0(\mod p-1)$}.\tag 21$$ If
$t\equiv 0 (\mod p-1)$, then $l_{p,q}(s_1, w^t)\equiv l_{p,q}(s_2,
w^t) (\mod p)$ for all $s_1, s_2 \in\Bbb Z_p ,$ $l_{p,q}(k,
w^t)\equiv l_{p,q}(k+p, w^t) (\mod p).$ It is easy to see that
$$\frac{1}{r+k-1}\binom{-r}{k}\binom{1-r-k}{j}=\frac{-1}{j+k}\binom{-r}{k+j-1}\binom{k+j}{j}
, \tag22$$ for all positive integers $r, j, k$ with $j, k \geq 0$,
$j+k>0, $ and $r\neq 1-k .$ Thus, we note that
$$\frac{1}{r+k-1}\binom{-r}{k}\binom{1-r-k}{j}=\frac{1}{r-1}\binom{-r+1}{k+j}\binom{k+j}{j}
.\tag 22-1$$ From (22) and (22-1), we derive
$$\frac{r}{r+k}\binom{-r-1}{k}\binom{-r-k}{j}=\binom{-r}{k+j}\binom{k+j}{j}.
\tag23$$ By using (13), we see that
$$\aligned
&\sum_{l=0}^{n-1}\frac{(-1)^{Fl+a}}{[Fl+a]_q^r}=\sum_{l=0}^{n-1}(-1)^l(-1)^a
\left([a]_q+q^a[F]_q[l]_{q^F}\right)^{-r}\\
&=-\sum_{s=0}^{\infty}[a]_q^{-r}\left(\frac{[F]_q}{[a]_q}\right)^sq^{as}(-1)^a\binom{-r}{s}
\frac{(-1)^n}{2}\sum_{l=0}^{s-1}\binom sl
q^{nFl}E_{l,q^F}[n]_{q^F}^{s-l}
\\
&-\sum_{s=0}^{\infty}[a]_q^{-r}\left(\frac{[F]_q}{[a]_q}\right)^sq^{as}(-1)^a\binom{-r}{s}
\frac{((-q^{Fs})^n-1)}{2}E_{s, q^F}.
\endaligned\tag24$$
For $s \in \Bbb Z_p ,$ we define the below $T$-Euler polynomials:
$$T_{n,q}(s,
a:F)=(-1)^a<a>^{-s}\sum_{k=0}^{\infty}\binom{-s}{k}[\frac{a}{F}]_{q^{F}}^{-k}q^{ak}
\left((-1)^nq^{nFk}-1\right)E_{k,q^F}.\tag25$$ Note that
$\lim_{q\rightarrow 1}T_{n,q}(s,a:F)=0,$ if $n$ is even positive
integer. From (23) and (24), we derive
$$\aligned
&\sum_{l=0}^{n-1}\frac{(-1)^{Fl+a}}{[Fl+a]_q^r}\\
&=-\sum_{s=0}^{\infty}\binom{-r}{s}
[a]_q^{-r}\left(\frac{[F]_q}{[a]_q}\right)^s \frac{
(-q^{s})^a(-1)^n}{2}\sum_{l=0}^{s-1}\binom{s}{l}q^{nFl}E_{l,
q^F}[n]_{q^F}^{s-l}\\
&-\frac{w^{-r}(a)}{2}T_{n,q}(r, a: F). \endaligned\tag26$$ First,
we evaluate the right side of Eq.(26) as follows:
$$\aligned
&\sum_{s=0}^{\infty}\binom{-r}{s}
[a]_q^{-r}\left(\frac{[F]_q}{[a]_q}\right)^s \frac{
(-q^{s})^a(-1)^n}{2}\sum_{l=0}^{s-1}\binom{s}{l}q^{nFl}E_{l,
q^F}[n]_{q^F}^{s-l}\\
 &=\sum_{k=0}^{\infty}\frac{r}{r+k}\binom{-r-1}{k}[a]_q^{-k-r}q^{ak}(-1)^n[Fn]_q^k
 \frac{(-1)^a}{2}\sum_{l=0}^{\infty}\binom{-r-k}{l}q^{al}
 \left(\frac{[F]_q}{[a]_q}\right)^l
 E_{l,q^F}.
\endaligned\tag27$$
It is easy to check that
$$q^{nFl}=\sum_{j=0}^l\binom lj
[nF]_q^j(q-1)^j=1+\sum_{j=1}^l\binom lj [nF]_q^j(q-1)^j. \tag 28$$
Let
$$K_{p,q}(s,
a:F)=\frac{(-1)^a}{2}<a>^{-s}\sum_{l=0}^{\infty}\binom{-s}{l}q^{al}\left(\frac{[F]_q}{[a]_q}\right)^l
E_{l,q^F}\sum_{j=1}^l\binom {l}{j}[nF]_q^j(q-1)^j.\tag29$$ Note
that $\lim_{q\rightarrow 1}K_{p,q}(s,a;F)=0.$
 For $F=p$,
$r\in\Bbb N$, we see that
$$2\sum_{a=1}^{p-1}\sum_{l=0}^{n-1}\frac{(-1)^{a+pl}}{[a+pl]_q^r}=2\sum_{\Sb
j=1 \\ (j, p)=1\endSb}^{np}\frac{(-1)^j}{[j]_q^r}.\tag 30$$ For
$s\in\Bbb Z_p$, we define $p$-adic analytically continued function
on $\Bbb Z_p$ as
$$\aligned
&K_{p,q}(s,\chi)=2\sum_{a=1}^{p-1}\chi(a)K_{p,q}(s,a:F),\\
&T_{p,q}(s,\chi)=2\sum_{a=1}^{p-1}\chi(a)T_{n,q}(s,a:F), \text{
where $k,n \geq 1$ .}
\endaligned\tag31$$
From (24)-(31), we derive
$$\aligned
&2\sum_{\Sb j=1 \\ (j,
p)=1\endSb}^{np}\frac{(-1)^j}{[j]_q^r}=-\sum_{k=0}^{\infty}\frac{r}{r+k}\binom{-r-1}{k}(-1)^n
[pn]_q^kl_{p,q}(r+k, w^{-r-k})\\
&-\sum_{k=0}^{\infty}\frac{r}{r+k}\binom{-r-1}{k}(-1)^n[pn]_q^kK_{p,q}(r+k,
w^{-r-k})-T_{p,q}(r, w^{-r}).
\endaligned$$
Therefore we obtain the following theorem: \proclaim{Theorem 5}
Let $p$ be an odd prime and let $n\geq 1$, and $r\geq 1$ be
integers. Then we have
$$\aligned
&2\sum_{\Sb j=1 \\ (j,
p)=1\endSb}^{np}\frac{(-1)^j}{[j]_q^r}=-\sum_{k=0}^{\infty}\frac{r}{r+k}\binom{-r-1}{k}(-1)^n
[pn]_q^kl_{p,q}(r+k, w^{-r-k})\\
&-\sum_{k=0}^{\infty}\frac{r}{r+k}\binom{-r-1}{k}(-1)^n[pn]_q^kK_{p,q}(r+k,
w^{-r-k})-T_{p,q}(r, w^{-r}).
\endaligned\tag32$$
\endproclaim
For $q=1$ in (32), we have
$$2\sum_{\Sb j=1 \\
(j,p)=1\endSb}^{np}\frac{(-1)^j}{j^r}=-\sum_{k=0}^{\infty}\frac{r}{k+r}\binom{-r-1}{k}(-1)^n(pn)^kl_p(r+k,
w^{-r-k}), $$ where $n$ is positive even integer.
 \proclaim{Remark } Let $p$ be an odd prime. Then we have
$$\sum_{j=1}^{p-1}\frac{(-1)^jq^j}{[j]_q}=\sum_{j=1}^{p-1}\frac{(-1)^j}{[j]_q}.$$
Proof. To prove Remark, it is sufficient to show that
$$\aligned
 \left(\sum_{j+1}^{p-1} \frac{(-1)^j}{[j]_q}
\right)^2&=\left(\sum_{j=1}^{p-1}\frac{(-1)^j}{[j]_q}\right)\left(\sum_{j=1}^{p-1}\frac{(-1)^j}{[j]_q}
-(1-q)\sum_{j=1}^{p-1}(-1)^j \right)\\
&=\left(\sum_{j=1}^{p-1}\frac{(-1)^j}{[j]_q}\right)\left(\sum_{j=1}^{p-1}(-1)^j\left(\frac{1}{[j]_q}-(1-q)\right)\right)\\
&=\left(\sum_{j=1}^{p-1}\frac{(-1)^j}{[j]_q}\right)\left(\sum_{j=1}^{p-1}\frac{(-1)^j
q^j }{[j]_q}\right).
\endaligned$$
\endproclaim


\Refs \ref \no 1 \by T. Kim \pages \paper A note on the
alternating  sums of powers of consecutive $q$-integers \yr  \vol
\jour (submitted)
\endref

\ref \key 2 \by T. Kim \pages 189-196 \paper Power series and
asymptotic series associated with the $q$-analogue of the two
variable $p$-adic $L$-function \yr 2005 \vol 12 \jour Russian J.
Math. Phys.\
\endref

\ref \key 3 \by T. Kim \pages 288-299 \paper $q$-Volkenborn
Integration \yr 2005 \vol 9 \jour Russ. J. Math. Phys.
\endref

\ref \key 4 \by T. Kim \pages 320-329 \paper On a $q$-analogue of
the $p$-adic log gamma functions \yr 1999 \vol 76 \jour J. Number
Theory
\endref

\ref \key 5 \by T. Kim \pages  \paper A note on the alternating
sums of powers of consecutive $q$-integers \yr 2006 \vol 1\jour
arXiv:math/0604227
\endref

\ref \key 6 \by T.Kim   \pages 15-18\paper Sums powers of
consecutive $q$-integers \yr 2004 \vol 9\jour Advan. Stud.
Contemp. Math.\endref

\ref\key 7 \by T. Kim \pages 179-187 \paper On $p$-adic
$q$-$L$-functions and sums of powers \yr 2002 \vol 252 \jour
Discrete Math.\endref

\ref\key 8 \by T. Kim \pages  \paper Multiple $p$-adic
$L$-functions\yr 2006 \vol 13 (2)\jour Russian J. Math.
Phys.\endref

\ref\key 9\by T. Kim \pages 13-17 \paper A note on $q$-Volkenborn
integration \yr 2005 \vol 8 \jour Proc. Jangjeon Math. Soc.\endref

\ref\key 10\by L. C. Washington \pages 50-61 \paper $p$-adic
$L$-functions and sums of powers \yr 1998 \vol 69 \jour J. Number
Theory
\endref

\ref\key 11\by L. C. Washington \pages  \book Introduction to
cyclotomic fields \yr 1982 \vol \publ Springer-Verlag(1'st Ed.)
\endref

\endRefs
\enddocument